\documentclass[12pt,reqno]{amsart} 

\usepackage{mathrsfs,amsfonts,amsthm,amsmath,amssymb} 
\usepackage[dvipsnames]{xcolor} 

\usepackage{enumitem}                
\usepackage{lineno}                 

\usepackage[numbers,sort&compress]{natbib}  

\usepackage{graphics}
\usepackage{graphicx}   
\graphicspath{{images/}} 
\usepackage{subfig}     
\usepackage{rotating}

\usepackage[hyperindex, pdfencoding = auto, psdextra, bookmarksdepth = 4]{hyperref}  
\hypersetup{
    bookmarksopen      = true,
    bookmarksopenlevel = 1,
    bookmarksnumbered  = true,
    pdftitle    = {},
    pdfauthor   = {},
    linktoc     = page,
    anchorcolor = green,
    colorlinks  = true,        
    linkcolor   = black,
    citecolor   = blue,
    filecolor   = blue,
    urlcolor    = magenta,
}

\usepackage[capitalise, noabbrev, nameinlink]{cleveref}    

\newtheorem{innercustomthm}{Theorem}
\newenvironment{customthm}[1]
    {\renewcommand\theinnercustomthm{#1}\innercustomthm}
    {\endinnercustomthm}

\theoremstyle{plain}

\newtheorem{theorem}{Theorem}[section]
\newtheorem{lemma}[theorem]{Lemma}

\newtheorem{proposition}[theorem]{Proposition}
\newtheorem{corollary}[theorem]{Corollary}
\newtheorem{question}[theorem]{Question}

\newtheorem{problem}[theorem]{Problem}

\theoremstyle{definition}   
\newtheorem{definition}[theorem]{Definition}
\newtheorem{remark}[theorem]{Remark}
\newtheorem{example}[theorem]{Example}

\makeatletter
    \def\subsection{\@startsection{subsection}{2}%
    \z@{.5\linespacing\@plus.7\linespacing}{.3\linespacing}%
    {\normalfont\bfseries}}
\makeatother 

\makeatletter                                         
    \newcommand{\LeftEqNo}{\let\veqno\@@leqno}        
\makeatother

\numberwithin{equation}{section}                      

\begin{document}

\

\vspace{-2cm}

\title[Strongly reducible operators in factors]{Strongly reducible operators in factors}

\author{Minghui Ma}
\address{Minghui Ma, School of Mathematical Sciences, Dalian University of Technology, Dalian, 116024, China}
\email{minghuima@dlut.edu.cn}

\author{Rui Shi}
\address{Rui Shi, School of Mathematical Sciences, Dalian University of Technology, Dalian, 116024, China}
\email{ruishi@dlut.edu.cn}

\author{Jianming Yang}
\address{Jianming Yang, School of Mathematical Sciences, Dalian University of Technology, Dalian, 116024, China}
\email{yangjianming@mail.dlut.edu.cn}

\author{Shanshan Yang}
\address{Shanshan Yang, School of Mathematical Sciences, Dalian University of Technology, Dalian, 116024, China}
\email{yss@mail.dlut.edu.cn}

\keywords{Strong reducibility, generalized kernel, algebraic operator, multi-cyclic operator, factor}

\makeatletter
\@namedef{subjclassname@2020}{\textup{2020} Mathematics Subject Classification}
\makeatother
\subjclass[2020]{Primary 46L10; Secondary  47C15}

\begin{abstract}
    Let $\mathcal M$ be a factor acting on a  separable complex Hilbert space. For every operator $T$ in $\mathcal M$ with dense generalized kernel and $T^2\ne 0$, we characterize the strong reducibility of $T$ relative to $\mathcal M$ based on a $3\times 3$ asymmetric matrix representation of $T$. As applications, we obtain criteria for the strong reducibility of algebraic operators in finite factors and type $\mathrm{III}$ factors. We also prove that every Hilbert space operator with an infinite-dimensional multicyclic direct summand is similar to an irreducible operator.
\end{abstract}

\maketitle

\section{Introduction}

\subsection{Strong reducibility and similarity orbits of operators}

The investigation of irreducible operators was initiated by Halmos in \cite{Halmos-68}, who established the uniform density of the set of irreducible operators on separable Hilbert spaces.
In 1970, Halmos \cite{Halmos-70} proposed ten problems in Hilbert space, the eighth of which asked whether the set of reducible operators is uniformly dense on each complex infinite-dimensional separable Hilbert space.
This problem was answered affirmatively by Voiculescu in his groundbreaking paper \cite{Voiculescu-76}. Whereas these results concern the topological properties of irreducible and reducible operators, Fillmore, Topping, and Radjavi, among others, studied irreducible operators from an algebraic perspective by investigating their sums \cite{FT69,Rad69}.

Note that irreducibility is invariant under unitary equivalence but not under similarity. Unlike unitary equivalence, similarity preserves the algebraic structure of the commutant but need not preserve orthogonal projections.  This distinction leads naturally to the question of which similarity orbits contain irreducible operators, or equivalently, which operators remain reducible under every similarity transformation. Similarity orbits have long been central to the structure and approximation theory of Hilbert space operators. The work of Apostol, Herrero, Voiculescu, and others connects approximation with spectral, Fredholm, and algebraic invariants \cite{Herrero-82,AFHV-84}. 

In 1972, Gilfeather \cite{Gilfeather-72} introduced a class of operators invariant under similarity, called \emph{strongly reducible} operators, and proved that every normal operator without eigenvalues is not strongly reducible. Recall that an operator in $\mathcal{B}(\mathcal{H})$ is strongly reducible if each operator in its similarity orbit is reducible. Later, Fong and Jiang \cite{FJ-94} characterized the strong reducibility of normal operators. Subsequently,  Hsin \cite{Hsin-00} obtained a complete characterization in finite dimensions.
Recently, the strong reducibility of normal operators in separable factors was treated in \cite{MSY-26}.
For further results on similarity and irreducibility, the reader is referred to \cite{JW-98-book,JW-06-book}.
This line of research naturally leads to the following question.
\begin{question} \label{ques-SR-B(H)}
    What is a necessary and sufficient condition for an operator in $\mathcal{B}(\mathcal{H})$ to be strongly reducible?
\end{question}
While normal operators admit a rich spectral theory, spectral data alone are insufficient to characterize strong reducibility for certain other classes of operators, notably nilpotent operators. Any complete answer to \Cref{ques-SR-B(H)} must nevertheless account for these classes. This necessitates new methods and structural invariants that capture features beyond the spectrum. In this paper, we investigate the strong reducibility of operators with dense generalized kernels, a class that includes strictly upper triangular operators, nilpotent operators, and suitable scalar translates of Cowen-Douglas operators.

\subsection{Operators with dense generalized kernels}

Throughout, let $\mathcal M$ be a factor acting on a separable Hilbert space $\mathcal H$, which is a von Neumann algebra with trivial center.
For an operator $T$ in $\mathcal M$, denote by $N(T)$ and $R(T)$ the \emph{null projection} and \emph{range projection} of $T$, i.e., projections onto $\ker T$ and $\overline{\operatorname{ran}T}$, respectively.
It is clear that $N(T) = I-R(T^*)$ and $N(T),R(T)\in W^*(T)\subseteq\mathcal M$, where $W^*(T)$ is the von Neumann algebra generated by $T$ and the identity $I$.
Recall that an operator in $\mathcal M$ is \emph{reducible} if it commutes with a nontrivial projection in $\mathcal M$.
Otherwise, it is \emph{irreducible} relative to $\mathcal M$.
An operator $T$ in a factor $\mathcal M$ is called \emph{strongly reducible} if $XTX^{-1}$ is reducible for every invertible operator $X$ in $\mathcal M$. An operator $T$ in $\mathcal M$  has \emph{a dense generalized kernel} if $\bigcup^{\infty}_{n=1} \operatorname{ker}T^n$ is dense in $\mathcal{H}$, equivalently, $\bigvee^{\infty}_{n=1} N(T^n) = I$. 

There are several related decomposition results.
By the Schur decomposition, every matrix in  $ M_n(\mathbb{C})$ is a sum of a normal matrix and a nilpotent one.
In this direction, Ringrose \cite{Ring-62} studied super-diagonal forms of compact operators on complex Banach spaces. An immediate consequence of his theorems is that every compact operator in $\mathcal{B}(\mathcal{H})$ is a sum of a normal operator and a quasinilpotent operator.  An analogous decomposition was established in \cite[Theorem 6]{DSZ15} for a diffuse finite von Neumann algebra $\mathcal{N}$ equipped with a faithful normal tracial state: every operator in  $\mathcal{N}$ is a sum of a normal operator and a s.o.t.-quasinilpotent operator. Through a different method, as a corollary of \cite[Theorem 1.3]{CFY24}, every operator in a type $\mathrm{II}_1$ factor is a sum of a normal operator and a nilpotent operator.


Note that nilpotent operators of order $2$ in a factor $\mathcal M$ with $\dim\mathcal M>4$ are already known to be strongly reducible by \cite[Lemma 4.3]{MSY-26}. Thus, it is sufficient to study operators $T$ satisfying $T^2 \neq 0$.
In this paper, we develop new methods to study the strong reducibility of operators with dense generalized kernels.
Throughout, when no confusion can arise, we assume that $T$ is an operator in $\mathcal M$ satisfying
\begin{equation}\label{equ assumption}
	T^2\ne 0\qquad\text{and}\qquad\bigvee_{n=1}^{\infty}N(T^n)=I.
\end{equation}
 By the invariant subspace structure, the strong reducibility of each operator satisfying \eqref{equ assumption} can be characterized with respect to a $3\times 3$ asymmetric matrix representation.

\subsection{A \texorpdfstring{$3\times 3$}{3x3} asymmetric matrix representation and main results}

To proceed, we introduce the matrix representation of $T$ given in \eqref{equ T-matrix}, using different decompositions of its domain and codomain.
This asymmetric matrix plays an important role in characterizing the strong reducibility of operators in separable factors. We prepare the following notation for the matrix in \eqref{equ T-matrix}.
Define projections 
\begin{equation}\label{equ EFG}
	E_T=N(T)-R\bigl(TN(T^2)\bigr),\qquad F_T=I-N(T^2),\qquad G_T=I-N(T),
\end{equation}
and operators
\begin{equation}\label{equ ABCD}
\begin{aligned}
	A_T &=E_TTF_T,\qquad C_T = R\bigl(TN(T^2)\bigr)T\bigl(N(T^2)-N(T)\bigr),\\
	B_T &=G_TTF_T,\qquad D_T = R\bigl(TN(T^2)\bigr)TF_T.
\end{aligned}
\end{equation}
With respect to the domain decomposition $I=N(T)+\bigl(N(T^2)-N(T)\bigr)+F_T$ and the codomain decomposition $I=E_T+R\bigl(TN(T^2)\bigr)+G_T$, the operator $T$ can be expressed in the matrix form
\begin{equation}\label{equ T-matrix}
	T=
	\begin{pmatrix}
		0 & 0 & A_T\\
		0 & C_T & D_T\\
		0 & 0 & B_T
	\end{pmatrix}.
\end{equation}
Details of the matrix in \eqref{equ T-matrix} are given in \Cref{sec T-matrix}.
Our first goal is to characterize the strong reducibility of $T$ in terms of a possibly unbounded operator $J_T$.

\begin{customthm}{\ref{thm SR}}
Let $\mathcal M$ be a separable factor and $T$ an operator in $\mathcal M$ satisfying \eqref{equ assumption} and of the form in \eqref{equ T-matrix}.
Let $\mathcal D_T=\{\xi\in\operatorname{ran}E_T\colon A_T^*\xi\in\operatorname{ran}B_T^*\}$ and
\begin{equation*}
	J_T=(B_T^*|_{\operatorname{ran}Q})^{-1}A_T^*\colon\mathcal D_T(\subseteq\operatorname{ran}P)\to\operatorname{ran}Q,
\end{equation*}
where $P$ is the projection onto $\overline{\mathcal D_T}$ and $Q=R(B_T)$.
The following statements are equivalent:
\begin{enumerate}  [label= $(\arabic*)$, ref= \arabic*, leftmargin=*] 
    \item[$(1)$] $T$ is strongly reducible in $\mathcal M$;
    \item[$(2)$] $R(A_T+YB_T)<E_T$ for every $Y\in E_T\mathcal MG_T$;
    \item[$(3)$] $B_T=KA_T+L$ for some relatively compact operator $K\in G_T\mathcal ME_T$ and finite operator $L\in G_T\mathcal MF_T$ satisfying that $R(L)\prec E_T-R(A_T)$;
    \item[$(4)$] $N(J_T)\succ N(J_T^*)$ and $\chi_{(0,c]}(|J_T|)$ is finite for every $c>0$.
\end{enumerate}
The kernels of $J_T$ and $J_T^*$ are taken in $\operatorname{ran}P$ and $\operatorname{ran}Q$, respectively.
\end{customthm}

A concrete example involving $J_T$ related to unilateral weighted shift operators in $\mathcal B(\mathcal H)$ can be found in \Cref{eg shift}.
By the above condition (3), $T$ is not strongly reducible unless $B_T$ is relatively compact in $\mathcal M$ (see \Cref{subsec factors}).
In particular, $T$ is not strongly reducible when $\mathcal M$ is of type $\mathrm{III}$ by \Cref{cor SR-III}.
As an application of \Cref{thm SR}, we have the following criterion when $A_T$ has closed range.

\begin{customthm}{\ref{thm SR-closed-range}}
Let $\mathcal M$ be a separable semifinite factor and $T$ an operator in $\mathcal M$ satisfying \eqref{equ assumption} and of the form in \eqref{equ T-matrix}.
If $B_T$ is relatively compact in $\mathcal M$ and $A_T$ has closed range, then the following statements are equivalent:
\begin{enumerate}  [label= $(\arabic*)$, ref= \arabic*, leftmargin=*] 
    \item[$(1)$] $T$ is strongly reducible in $\mathcal M$;
    \item[$(2)$] $N(A_T)\wedge F_T\prec N(A_T^*)\wedge E_T$.
\end{enumerate}
\end{customthm}

If $\mathcal M=\mathcal B(\mathcal H)$, then the above condition (2) means that $A_T$ is a semi-Fredholm operator with $\operatorname{ind}A_T<0$ (allowing $-\infty$) as an operator in $\mathcal B(F_T\mathcal H,E_T\mathcal H)$.
If $\mathcal M$ is a type $\mathrm{II}_\infty$ factor, a generalized Fredholm theory was established in \cite{Breuer-68}.
In general, we have no simple criterion when $A_T$ does not have closed range.
In the following theorem, we consider the special case that $B_T$ is a finite operator.
Note that
\begin{equation*}
	R(B_T)\sim R(B_T^*)=F_T=R((T^*)^2)\sim R(T^2).
\end{equation*}

\begin{customthm}{\ref{thm SR-B-finite}}
Let $\mathcal M$ be a separable semifinite factor and $T$ an operator in $\mathcal M$ satisfying condition \eqref{equ assumption}.
If $T^2$ is a finite operator, then the following statements are equivalent:
\begin{enumerate} [label= $(\arabic*)$, ref= \arabic*, leftmargin=*] 
    \item[$(1)$] $T$ is strongly reducible in $\mathcal M$;
    \item[$(2)$] $R(T^2)\prec N(T)-R(TN(T^2))$; equivalently, $F_T\prec E_T$.
\end{enumerate}
\end{customthm}

In the following two theorems, we provide a criterion for the strong reducibility of algebraic operators in finite factors and type $\mathrm{III}$ factors.
Since every operator in the matrix algebra $M_n(\mathbb{C})$ is algebraic, the following result can be viewed as a generalization of Hsin's theorem \cite{Hsin-00}.

\begin{customthm}{\ref{thm Hsin}}
Let $(\mathcal M,\tau)$ be a separable finite factor with $\dim\mathcal M>4$ and $T$ an algebraic operator in $\mathcal M$.
Then $T$ is not strongly reducible if and only if
\begin{enumerate} [label= $(\arabic*)$, ref= \arabic*, leftmargin=*] 
    \item[$(1)$] $\tau(R(\lambda I-T))\geqslant\frac{1}{2}$ for every scalar $\lambda$ in $\mathbb{C}$;
    \item[$(2)$] the set $\{I,T,T^2\}$ is linearly independent over $\mathbb{C}$.
\end{enumerate}
\end{customthm}

Due to \Cref{cor SR-III}, we can drop the trace condition in \Cref{thm Hsin} for algebraic operators in type $\mathrm{III}$ factors.

\begin{customthm}{\ref{thm SR-III}}
Let $\mathcal M$ be a separable type $\mathrm{III}$ factor and $T$ an algebraic operator in $\mathcal M$.
Then $T$ is not strongly reducible if and only if the set $\{I,T,T^2\}$ is linearly independent over $\mathbb{C}$.
\end{customthm}

It is natural to ask the following problem.

\begin{problem}
Do the characterizations of strong reducibility in \Cref{thm Hsin} and \Cref{thm SR-III} hold for general operators in factors of type $\mathrm{II}_1$ and type $\mathrm{III}$?
\end{problem}

Recall that an operator $S\in\mathcal{B}(\mathcal{H})$ is said to be \emph{multi-cyclic} if there exists a \emph{generating set} $\{\xi_j\}_{j=1}^m$ of finitely many vectors in $\mathcal{H}$ such that
\begin{equation*}
	\mathcal{H}=\overline{\mathrm{span}}\{S^k\xi_j\colon k\geqslant 0,~1\leqslant j\leqslant m\}.
\end{equation*}
We say that $S$ is \emph{cyclic} if $m=1$.
Clearly, the set of multi-cyclic operators is invariant under similarity and the unilateral shift operator is cyclic.
The following theorem asserts that every operator with an infinite-dimensional multi-cyclic direct summand is similar to an irreducible operator.

\begin{customthm}{\ref{thm multi-cyclic}}
Let $S$ be a multi-cyclic operator acting on an infinite-dimensional Hilbert space $\mathcal{H}$.
Then $T\oplus S$ is not strongly reducible for every operator $T$ acting on a separable Hilbert space.
In particular, $S$ is not strongly reducible.
\end{customthm}

This paper is organized as follows.
In \Cref{sec 2}, we develop the matrix representation and establish a range criterion for strong reducibility.
In \Cref{sec SR}, we analyze this criterion using closed densely defined operators and prove \Cref{thm SR}.
\Cref{sec application} presents further criteria and applications to Cowen--Douglas operators and Jordan models.
In \Cref{sec algebric}, we characterize the strong reducibility of algebraic operators in finite factors and type $\mathrm{III}$ factors.
Finally, in \Cref{sec multi-cyclic}, we prove that every operator with an infinite-dimensional multi-cyclic direct summand is similar to an irreducible operator.

\section{Preliminaries}
\label{sec 2}
\subsection{Factors}
\label{subsec factors}

Let $\mathcal H$ be a Hilbert space and $\mathcal B(\mathcal H)$ the algebra of all bounded linear operators acting on $\mathcal H$.
A \emph{von Neumann algebra} is a weak-operator closed self-adjoint unital subalgebra of $\mathcal B(\mathcal H)$ and  a \emph{factor} is a von Neumann algebra with trivial center.
By the type decomposition theorem \cite[Theorem 6.5.2]{KR2}, factors are classified into type $\mathrm{I}_n$, type $\mathrm{I}_\infty$, type $\mathrm{II}_1$, type $\mathrm{II}_\infty$, and type $\mathrm{III}$.

Throughout, let $\mathcal M$ denote a factor acting on a separable Hilbert space $\mathcal H$.
For two projections $P,Q$ in $\mathcal M$, we say that $P$ is \emph{(Murray-von Neumann) equivalent} to $Q$ (denoted by $P\sim Q$) if $P=V^*V$ and $Q=VV^*$ for some partial isometry $V$ in $\mathcal M$, $P$ is \emph{weaker} than $Q$ (denoted by $P\precsim Q$ or $Q\succsim P$) if $P$ is equivalent to a subprojection of $Q$, and $P$ is \emph{strictly weaker} than $Q$ (denoted by $P\prec Q$ or $Q\succ P$) if $P\precsim Q$ and $P\nsim Q$.
A projection $P$ is \emph{finite} if every proper subprojection of $P$ is strictly weaker than $P$.
Otherwise, $P$ is said to be \emph{infinite}.
Note that $P\prec Q$ forces $P=0$ in every separable type $\mathrm{III}$ factor as all nonzero projections are equivalent.

Recall that the set $\mathcal F_{\mathcal M}$ of all \emph{finite operators} in $\mathcal M$ is the ideal generated by all finite projections in $\mathcal M$, and the set $\mathcal K_{\mathcal M}$ of \emph{relatively compact operators} is the norm closure of $\mathcal F_{\mathcal M}$.
In $\mathcal B(\mathcal H)$, finite operators are finite-rank operators and relatively compact operators are compact operators, respectively.
Note that $\mathcal K_{\mathcal M}=\{0\}$ if $\mathcal M$ is of type $\mathrm{III}$.

\subsection{Matrix representation}
\label{sec T-matrix}

In this subsection, we explain the matrix form \eqref{equ T-matrix} of $T$ in detail.
For every operator $T$ in $\mathcal M$, consider the following chain of projections
\begin{equation*}
	N(T)\leqslant N(T^2)\leqslant N(T^3)\leqslant\cdots\leqslant\bigvee_{n=1}^{\infty}N(T^n)\leqslant I,
\end{equation*}
where the union $\bigvee_{n=1}^{\infty}N(T^n)$ is the projection onto the closure of the \emph{generalized kernel} $\bigcup_{n=1}^{\infty}\ker T^n$.
Under the assumption of \eqref{equ assumption}, we define a sequence of projections with sum $I$ as follows
\begin{equation}\label{equ kernel-flag}
	P_1=N(T)\qquad\text{and}\qquad P_j=N(T^j)-N(T^{j-1})\qquad\text{for}~j\geqslant 2.
\end{equation}
Note that $P_j=0$ for all $j\geqslant m+1$ if $T$ is a nilpotent operator of order $m$.
With respect to the decomposition $I=\sum_{j=1}^{\infty}P_j$, $T$ has a strictly upper triangular matrix representation of the form
\begin{equation*}
	T=
	\begin{pmatrix}
		0 & T_{12} & T_{13} & \cdots\\
		& 0 & T_{23} & \cdots\\
		& & 0 & \cdots\\
		& & & \ddots
	\end{pmatrix}
	\begin{matrix}
		\operatorname{ran}P_1\\
		\operatorname{ran}P_2\\
		\operatorname{ran}P_3\\
		\vdots\\
	\end{matrix},
\end{equation*}
where $T_{ij}=P_iTP_j\in P_i\mathcal MP_j$ for all $i<j$, and $R(T_{j,j+1}^*)=P_{j+1}$ for all $j\geqslant 1$.
The projections defined in \eqref{equ EFG} can be written as
\begin{equation*}
	E_T=P_1-R(T_{12}),\qquad F_T=\sum_{j=3}^{\infty}P_j,\qquad G_T=\sum_{j=2}^{\infty}P_j,
\end{equation*}
and the corner operators defined in \eqref{equ ABCD} become
\begin{equation*}
	A_T=E_TTF_T,\qquad B_T=G_TTF_T,\qquad C_T=R(T_{12})TP_2,\qquad D_T=R(T_{12})TF_T.
\end{equation*}
Moreover, we have
\begin{equation*}
	R(B_T^*)=F_T,\qquad R(C_T^*)=P_2,\qquad R(C_T)=R(T_{12})=R\bigl(TN(T^2)\bigr).
\end{equation*}
We provide a formula for the common kernel of $T$ and $T^*$ as follows.

\begin{lemma}\label{lem common-kernel}
With the notation of \eqref{equ EFG} and \eqref{equ ABCD}, we have 
\begin{equation*}
	E_T-R(A_T)=N(T)\wedge N(T^*).
\end{equation*}
\end{lemma}

\begin{proof}
Since $E_TTN(T^2)=0$, we have $A_T=E_TTF_T=E_TT$, i.e., $A_T^*=T^*E_T$.
It follows that
\begin{equation*}
	E_T-R(A_T)=E_T\wedge N(A_T^*)=E_T\wedge N(T^*).
\end{equation*}
Since $N(T^*)=I-R(T)\leqslant I-R\bigl(TN(T^2)\bigr)$, we have
\begin{equation*}
	N(T)\wedge N(T^*)\leqslant N(T)-R\bigl(TN(T^2)\bigr)=E_T.
\end{equation*}
Thus, $N(T)\wedge N(T^*)\leqslant E_T\wedge N(T^*)$, and hence $E_T\wedge N(T^*)=N(T)\wedge N(T^*)$.
This completes the proof.
\end{proof}

\subsection{Criteria for reducibility and irreducibility}
\label{subsec criteria}

In this subsection, let $T$ be an operator in $\mathcal M$ satisfying condition \eqref{equ assumption}.
In the following two lemmas, we provide criteria for the reducibility and irreducibility of $T$ of the form \eqref{equ T-matrix}.

\begin{lemma}\label{lem R-criterion}
With the notation of \eqref{equ T-matrix}, $T$ is reducible if $R(A_T)<E_T$.
\end{lemma}

\begin{proof}
Let $Q=E_T-R(A_T)$.
Since $QT=0=TQ$ by \Cref{lem common-kernel}, $Q$ is a nontrivial projection commuting with $T$.
This completes the proof.
\end{proof}

\begin{lemma}\label{lem IR-criterion}
With the notation of \eqref{equ T-matrix}, $T$ is irreducible if
\begin{equation*}
	R(A_T)=E_T\qquad\text{and}\qquad\{C_TC_T^*,D_TD_T^*\}'\cap Q_0\mathcal MQ_0=\mathbb{C}Q_0,
\end{equation*}
where $Q_0=R(C_T)$.
\end{lemma}

\begin{proof}
Let $Q$ be a projection commuting with $T$.
Then $Q$ commutes with $Q_0$, $E_T$, and $P_j$ for all $j\geqslant 1$.
Since $C_TC_T^*,D_TD_T^*\in W^*(T)$, we have $Q_0Q=0$ or $Q_0Q=Q_0$ by assumption.
Replacing $Q$ by $I-Q$ if necessary, we assume that $Q_0Q=0$.
Since
\begin{equation*}
	R(C_T^*)=P_2,\qquad C_TQ=Q_0TP_2Q=Q_0QTP_2=0,
\end{equation*}
we have $P_2Q=0$.
For every $j\geqslant 2$, we have
\begin{equation*}
	R(T_{j,j+1}^*)=P_{j+1},\qquad T_{j,j+1}Q=P_jTP_{j+1}Q=P_jQTP_{j+1}.
\end{equation*}
It follows that $P_{j+1}Q=0$ whenever $P_jQ=0$.
Thus, $P_jQ=0$ for all $j\geqslant 2$ by induction.
In particular, $F_TQ=QF_T=0$ and hence
\begin{equation*}
	QA_T=QE_TTF_T=E_TTQF_T=0.
\end{equation*}
Thus, $QE_T=0$ by the assumption $R(A_T)=E_T$.
Therefore, $Q=0$ and hence $T$ is irreducible in $\mathcal M$.
\end{proof}

\subsection{A range condition}
\label{subsec range}

In this subsection, let $T$ be an operator in $\mathcal M$ satisfying condition \eqref{equ assumption}.
The following proposition gives a range criterion involving $A_T$ and $B_T$.
The implications $(2)\Rightarrow(1)$ and $(1)\Rightarrow(2)$ are proved in \Cref{cor SR-criterion} and \Cref{cor SR-criterion-converse}, respectively.

\begin{proposition}\label{prop SR-criterion}
With the notation of \eqref{equ T-matrix}, the following conditions are equivalent:
\begin{enumerate}
\item[$(1)$] $T$ is strongly reducible in $\mathcal M$;
\item[$(2)$] $R(A_T+YB_T)<E_T$ for every $Y\in E_T\mathcal MG_T$.
\end{enumerate}
\end{proposition}

Since the projection $R(A_T)$ is closely related to the reducibility and irreducibility of $T$ by \Cref{lem R-criterion,lem IR-criterion}, we need to describe how $A_T$ transforms under similarity.
To facilitate the argument, we introduce the following definition.

\begin{definition}\label{def corner-invertible}
Let $E$ and $F$ be projections in $\mathcal M$.
An operator $X$ in $E\mathcal MF$ is said to be \emph{corner-invertible} if there exists an operator $Y$ in $F\mathcal ME$, called the \emph{corner-inverse} of $X$, such that $XY=E$ and $YX=F$.
\end{definition}

\begin{lemma}\label{lem A-transform}
Let $X$ be an invertible operator in $\mathcal M$ and $S=XTX^{-1}$.
With the notation of \eqref{equ T-matrix}, there are corner-invertible operators $Y_1\in E_S\mathcal ME_T$, $Y_2\in F_T\mathcal MF_S$, and an operator $Y\in E_T\mathcal MG_T$ such that
\begin{equation*}
	A_S=Y_1(A_T+YB_T)Y_2.
\end{equation*}
\end{lemma}

\begin{proof}
It is clear that $\ker S^j=X\ker T^j$ for all $j\geqslant 1$.
As in \eqref{equ kernel-flag}, let
\begin{equation*}
	Q_1=N(S)\qquad\text{and}\qquad Q_j=N(S^j)-N(S^{j-1})\qquad\text{for}~j\geqslant 2.
\end{equation*}
Then $Q_iXP_j=0$ and $P_iX^{-1}Q_j=0$ whenever $i>j$.
It is clear that
\begin{equation*}
	(P_jX^{-1}Q_j)(Q_jXP_j)=P_j\qquad\text{and}\qquad(Q_jXP_j)(P_jX^{-1}Q_j)=Q_j
\end{equation*}
for all $j\geqslant 1$.
Thus, with respect to the domain decomposition $I=\sum_{j=1}^{\infty}P_j$ and the codomain decomposition $I=\sum_{j=1}^{\infty}Q_j$, the operator $X$ is upper triangular with each diagonal element corner-invertible.
For simplicity, let
\begin{equation*}
	J=E_SXE_T,\qquad K=E_SXG_T,\qquad L=F_SXF_T.
\end{equation*}
A direct computation shows that $\overline{\operatorname{ran}S_{12}}=X\overline{\operatorname{ran}T_{12}}$.
Thus, $J$ is corner-invertible in $E_S\mathcal ME_T$ with corner-inverse $E_TX^{-1}E_S$, and $L$ is corner-invertible in $F_S\mathcal MF_T$ with corner inverse $F_TX^{-1}F_S$.
By $A_S=E_SSF_S=E_SS$, $SX=XT$, and $E_SXR(T_{12})=0$, we have
\begin{align*}
	A_SL &=E_SSXF_T=E_SXTF_T=E_SX\bigl(I-R(T_{12})\bigr)TF_T\\
	&=E_SX(A_T+B_T)=JA_T+KB_T.
\end{align*}
We complete the proof by taking $Y_1=J\in E_S\mathcal ME_T$, $Y_2=L^{-1}\in F_T\mathcal MF_S$, and $Y=J^{-1}K\in E_T\mathcal MG_T$, where we denote by $A^{-1}$ the corner-inverse of an operator $A$ when no confusion can arise.
\end{proof}

The following corollary follows immediately from \Cref{lem R-criterion} and \Cref{lem A-transform}.

\begin{corollary}\label{cor SR-criterion}
With the notation of \eqref{equ T-matrix}, $T$ is strongly reducible in $\mathcal M$ if
\begin{equation*}
	R(A_T+YB_T)<E_T\qquad\text{for every}~Y\in E_T\mathcal MG_T.
\end{equation*}
\end{corollary}

Next we present a technical lemma that will be used in the proof of \Cref{cor SR-criterion-converse}.

\begin{lemma}\label{lem CC*-DD*}
Let $F$, $G$, $P$ and $Q$ be nonzero projections in $\mathcal M$.
Suppose that the operators $B\in G\mathcal MF$, $C\in Q\mathcal MP$, and $D\in Q\mathcal MF$ satisfy 
\begin{equation*}
	B\ne 0,\qquad R(C^*)=P,\qquad R(C)=Q.
\end{equation*}
Then there exist a positive invertible operator $H\in P\mathcal MP$ and an operator $K\in Q\mathcal MG$ such that
\begin{equation*}
	\{CHC^*,(D+KB)(D+KB)^*\}'\cap Q\mathcal MQ=\mathbb{C}Q.
\end{equation*}
\end{lemma}

\begin{proof}
Multiplying by a positive scalar, we may assume that $\|C\|=1$.
Let $C=AU$ be the polar decomposition, where $A=(CC^*)^{1/2}$, $U^*U=P$, and $UU^*=Q$.
For every $n\geqslant0$, let
\begin{equation*}
	E_n=\chi_{(2^{-n-1},2^{-n}]}(A).
\end{equation*}
Then $Q=\sum_{n=0}^{\infty}E_n$ and $E_0\ne0$.
Since $\mathcal M$ is a factor, we can choose a nonzero partial isometry $V_0\in E_0\mathcal MR(B)$ and, for every $n\geqslant1$ with $E_n\ne0$, a nonzero partial isometry $V_n\in E_n\mathcal MV_0^*V_0$, requiring that
\begin{equation*}
	V_nV_n^*<E_n
	\qquad\text{for every}~n\geqslant0~
	\text{with}~E_n\mathcal ME_n\ne\mathbb CE_n.
\end{equation*}
Then $V_nB\ne0$ for every nonzero $E_n$.
It follows from $V_n=V_nV_0^*V_0$ that $V_0BB^*V_n^*\ne0$.
Let
\begin{equation*}
	V=\sum_{E_n\ne0}2^{-n}V_n,\qquad
	X_n(t)=E_nD+t2^{-n}V_nB\quad(t\in\mathbb R).
\end{equation*}
The series defining $V$ converges in norm and the quadratic coefficient of $X_0(t)X_n(t)^*$ is nonzero.
Whenever $E_n\mathcal ME_n\ne\mathbb CE_n$, the quadratic coefficient of $X_n(t)X_n(t)^*$ is nonzero and supported on $V_nV_n^*<E_n$, and hence is nonscalar.
Thus, each condition
\begin{equation*}
	X_0(t)X_n(t)^*=0
	\qquad\text{or}\qquad
	X_n(t)X_n(t)^*\in\mathbb RE_n
\end{equation*}
holds for at most two real values of $t$, with the second condition considered only when $E_n\mathcal ME_n\ne\mathbb CE_n$.
Choose $0<t<1$ outside the countable union of these exceptional sets.
Let $K=tV$ and $X_n=E_n(D+KB)$.
Then
\begin{equation}\label{equ 0-n}
	X_0X_n^*\ne0\qquad\text{whenever}~E_n\ne0,
\end{equation}
and $Y_n=X_nX_n^*$ is nonscalar whenever $E_n\mathcal ME_n\ne\mathbb CE_n$.
By \cite[Corollary 3.12]{MSY-26}, and trivially in the scalar case, we can choose $Z_n\in(E_n\mathcal ME_n)_{sa}$ such that
\begin{equation}\label{equ Z}
	\sigma_{E_n\mathcal ME_n}(Z_n)\subseteq(4^{-n-1},4^{-n}),
	\qquad
	\{Y_n,Z_n\}'\cap E_n\mathcal ME_n=\mathbb CE_n.
\end{equation}
Let
\begin{equation*}
	Z=\sum_{E_n\ne0}Z_n,\qquad
	H_0=\sum_{E_n\ne0}(AE_n)^{-1}Z_n(AE_n)^{-1},
\end{equation*}
where the inverses are taken in $E_n\mathcal ME_n$ and the sums converge strongly.
Since
\begin{equation*}
	\frac14E_n\leqslant(AE_n)^{-1}Z_n(AE_n)^{-1}\leqslant4E_n,
\end{equation*}
$H_0$ is positive and invertible in $Q\mathcal MQ$.
Let $H=U^*H_0U$.
Then $H$ is positive and invertible in $P\mathcal MP$, and
\begin{equation*}
	CHC^*=AH_0A=Z.
\end{equation*}
Let $W=(D+KB)(D+KB)^*$ and let $R\in Q\mathcal MQ$ be a projection commuting with $Z$ and $W$.
The disjoint spectral intervals in \eqref{equ Z} imply that $R$ commutes with every $E_n$.
Since $E_nR$ commutes with $Z_n$ and $Y_n=E_nWE_n$, \eqref{equ Z} gives
\begin{equation*}
	R=\sum_{E_n\ne0}\varepsilon_nE_n,\qquad \varepsilon_n\in\{0,1\}.
\end{equation*}
Commuting with $W$ yields $(\varepsilon_0-\varepsilon_n)X_0X_n^*=0$.
Thus, all $\varepsilon_n$ are equal by \eqref{equ 0-n}, and hence $R=0$ or $R=Q$.
This completes the proof.
\end{proof}

Now we are ready to prove the converse of \Cref{cor SR-criterion}.

\begin{corollary}\label{cor SR-criterion-converse}
With the notation of \eqref{equ T-matrix}, if $T$ is strongly reducible in $\mathcal M$, then
\begin{equation*}
	R(A_T+YB_T)<E_T\qquad\text{for every}~Y\in E_T\mathcal MG_T.
\end{equation*}
\end{corollary}

\begin{proof}
Suppose that $R(A_T+YB_T)=E_T$ for some $Y\in E_T\mathcal MG_T$.
We claim that $T$ is not strongly reducible.
Since $Y^2=TY=0$, the operator $I+Y$ is invertible in $\mathcal M$ with inverse $I-Y$.
Let $S_0=(I+Y)T(I+Y)^{-1}$.
Combining this with \eqref{equ T-matrix}, we have 
\begin{equation*}
	S_0=
	\begin{pmatrix}
		0 & 0 & A_T+YB_T\\
		0 & C_T & D_T\\
		0 & 0 & B_T
	\end{pmatrix}.
\end{equation*}
It follows that $R(A_{S_0})=E_T=E_{S_0}$.
Up to similarity, we may assume that
\begin{equation*}
	R(A_T)=E_T.
\end{equation*}
By \Cref{lem CC*-DD*}, there exist a positive invertible operator $H\in P_2\mathcal MP_2$ and an operator $K\in Q_0\mathcal MG_T$ such that
\begin{equation*}
	\{C_THC_T^*,(D_T+KB_T)(D_T+KB_T)^*\}'\cap Q_0\mathcal MQ_0=\mathbb{C}Q_0,
\end{equation*}
where $Q_0=R(C_T)$.
Note that $K^2=TK=0$.
Let $S=XTX^{-1}$, where $X$ is an invertible operator given by
\begin{equation*}
	X=(P_1+H^{-1/2}+F_T)(I+K),\qquad
	X^{-1}=(I-K)(P_1+H^{1/2}+F_T).
\end{equation*}
As in the proof of \Cref{lem A-transform}, let
\begin{equation*}
	Q_1=N(S),\qquad Q_j=N(S^j)-N(S^{j-1})\quad\text{for}~j\geqslant 2.
\end{equation*}
Then $Q_1=P_1$, $Q_2=P_2$, $R(C_S)=R(C_T)=Q_0$, and hence $E_S=E_T$, $F_S=F_T$.
It follows that $A_S=A_T$, $C_S=C_TH^{1/2}$, and $D_S=D_T+KB_T$.
In particular, we have
\begin{equation*}
	R(A_S)=E_S,\qquad\{C_SC_S^*,D_SD_S^*\}'\cap Q_0\mathcal MQ_0=\mathbb{C}Q_0.
\end{equation*}
By \Cref{lem IR-criterion}, $S$ is irreducible in $\mathcal M$.
Therefore, $T$ is not strongly reducible.
\end{proof}

\section{Strongly reducible operators}
\label{sec SR}

In this section, we prove \Cref{thm SR}, which gives a complete characterization of the strong reducibility of operators satisfying \eqref{equ assumption}.
To this end, we analyze the range condition in \Cref{prop SR-criterion}.

Suppose that $E,F,G$ are projections in $\mathcal M$, and $A\in E\mathcal MF$, $B\in G\mathcal MF$.
We investigate the following range condition
\begin{equation}\label{equ property-UDR}
	R(A+YB)<E\qquad\text{for every}~Y\in E\mathcal MG.
\end{equation}
We shall see in \Cref{lem J-AB} that the range condition \eqref{equ property-UDR} is naturally related to a densely defined closed operator $J_{A,B}$ affiliated with $\mathcal M$.

Let $Q=R(B)\leqslant G$ be the range projection of $B$.
Then $\operatorname{ran}Q=\mathcal H\ominus\ker B^*$ and
\begin{equation*}
	B^*|_{\operatorname{ran}Q}\colon\operatorname{ran}Q\to\overline{\operatorname{ran}B^*}
\end{equation*}
is an injective operator with dense range.
Let $\mathcal D_{A,B}=\{\xi\in\operatorname{ran}E\colon A^*\xi\in\operatorname{ran}B^*\}$ be a linear subspace of $\operatorname{ran}E$, $P$ the projection onto $\overline{\mathcal D_{A,B}}$, and
\begin{equation}\label{equ J-AB}
	J_{A,B}=(B^*|_{\operatorname{ran}Q})^{-1}A^*\colon\mathcal D_{A,B}(\subseteq\operatorname{ran}P)\to\operatorname{ran}Q.
\end{equation}
It is clear that
\begin{equation*}
	\operatorname{Graph}(J_{A,B})=\{(\xi,\eta)\in\operatorname{ran}P\oplus\operatorname{ran}Q\colon A^*\xi=B^*\eta\}.
\end{equation*}
Therefore, $P\in E\mathcal ME$, $Q\in G\mathcal MG$, and $J_{A,B}$ is a densely defined closed operator affiliated with $Q\mathcal MP$, i.e., $QJ_{A,B}P$ is affiliated with $\mathcal M$.
Note that $N(J_{A,B})$ is understood as the projection onto $\ker J_{A,B}$.
In particular, $N(J_{A,B})\in P\mathcal MP$.

\begin{lemma}\label{lem J-AB}
With the notation of \eqref{equ J-AB}, the following statements are equivalent:
\begin{enumerate}
\item[$(1)$] $R(A+YB)<E$ for every $Y\in E\mathcal MG$;
\item[$(2)$] $N(J_{A,B}+X)\ne 0$ for every $X\in Q\mathcal MP$.
\end{enumerate}
\end{lemma}

\begin{proof}
Let $Y\in E\mathcal MG$.
If $\xi\in\operatorname{ran}E$ satisfies $(A+YB)^*\xi=0$, then
\begin{equation*}
	A^*\xi=-B^*Y^*\xi\in\operatorname{ran}B^*,
\end{equation*}
and hence $\xi\in\mathcal D_{A,B}$.
For every $\xi\in\mathcal D_{A,B}$, we have
\begin{equation*}
	(A+YB)^*\xi
	=A^*\xi+B^*Y^*\xi
	=B^*(J_{A,B}+QY^*P)\xi.
\end{equation*}
Since $B^*|_{\operatorname{ran}Q}$ is injective, it follows that
\begin{equation*}
	E-R(A+YB)=N(J_{A,B}+QY^*P),
\end{equation*}
where the kernel on the right is taken in $\operatorname{ran}P$.
Finally, every $X\in Q\mathcal MP$ is of the form $QY^*P$ for some $Y\in E\mathcal MG$, by taking $Y=X^*$.
This proves the equivalence.
\end{proof}

We study condition $(2)$ in \Cref{lem J-AB} for a general densely defined closed operator $J$ affiliated with $Q\mathcal MP$, not necessarily defined by the operators $A$ and $B$ as in \eqref{equ J-AB}.
The kernels of $J$ and $J^*$ are taken inside $\operatorname{ran}P$ and $\operatorname{ran}Q$, i.e., $N(J)\leqslant P$ and $N(J^*)\leqslant Q$, respectively.

\begin{proposition}\label{prop UK}
Let $P,Q$ be projections in $\mathcal M$ and $J\colon\mathcal D(\subseteq\operatorname{ran}P)\to\operatorname{ran}Q$ a closed densely defined operator affiliated with $Q\mathcal MP$.
The following statements are equivalent:
\begin{enumerate}
\item[$(1)$] $N(J+X)\ne 0$ for every $X\in Q\mathcal MP$;
\item[$(2)$] $N(J)\succ N(J^*)$ and $\chi_{(0,c]}(|J|)$ is finite for every $c>0$.
\end{enumerate}
\end{proposition}

\begin{proof}
Let $J=U|J|$ be the polar decomposition, $E=N(J)$, and $F=N(J^*)$.

$(1)\Rightarrow(2)$.
If $E\precsim F$, then there exists a partial isometry $W\in F\mathcal ME$ such that $W^*W=E$ and $WW^*\leqslant F$.
Since the ranges of $J$ and $W$ are orthogonal and $W$ is isometric on $\operatorname{ran}E$, $J+W$ is injective, contradicting $(1)$.
Thus, $E\succ F$.

If $G=\chi_{(0,c]}(|J|)$ is infinite for some $c>0$, then both $E+G$ and $F+UGU^*$ are infinite projections.
Thus, there exists a partial isometry $V$ in $\mathcal M$ such that
\begin{equation*}
	V^*V=E+G\quad\text{and}\quad VV^*=F+UGU^*.
\end{equation*}
Since $JG$ is bounded, we can take $X=V-JG\in Q\mathcal MP$.
For every $\xi\in\mathcal D$, we have
\begin{equation*}
	(J+X)\xi=J(P-E-G)\xi+V(E+G)\xi.
\end{equation*}
Since $J(P-E-G)\xi\in\operatorname{ran}(F+UGU^*)^\perp$ and $V(E+G)\xi\in\operatorname{ran}(F+UGU^*)$, we see that $N(J+X)=0$.
Thus, $(1)$ implies that $\chi_{(0,c]}(|J|)$ is finite for every $c>0$.

$(2)\Rightarrow(1)$.
Let $X\in Q\mathcal MP$, $c>\|X\|$, and $G=\chi_{(0,c]}(|J|)$.
For simplicity, let
\begin{equation*}
	P_1=E+G,\quad P_2=P-P_1,\quad Q_1=F+UGU^*,\quad Q_2=Q-Q_1.
\end{equation*}
Since $E\succ F$ and $G\sim UGU^*$ is finite, we have $P_1\succ Q_1$.
Note that the restriction
\begin{equation*}
	J|_{\mathcal D\cap\operatorname{ran}P_2}\colon\mathcal D\cap\operatorname{ran}P_2\to\operatorname{ran}Q_2
\end{equation*}
has a bounded inverse $L\in P_2\mathcal MQ_2$ satisfying
\begin{equation*}
	\operatorname{ran}L\subseteq\mathcal D,\qquad
	JL=Q_2,\qquad \|L\|\leqslant c^{-1}.
\end{equation*}
Since $\|LX\|<1$, the operator $P+LX$ is invertible in $P\mathcal MP$, whose inverse is denoted by $Z$.
For every $\xi\in\operatorname{ran}P_1\subseteq\mathcal D$, we have $Z\xi=\xi-LXZ\xi\in\mathcal D$.
Since $JP_1$ is bounded and $JL=Q-Q_1$, it follows that
\begin{equation*}
	(J+X)ZP_1=JP_1+Q_1XZP_1\in Q_1\mathcal MP_1.
\end{equation*}
Since $P_1\succ Q_1$, we can choose $0\ne\xi\in\operatorname{ran}P_1$ such that $(J+X)Z\xi=0$.
Since $Z$ is invertible, $Z\xi\in\mathcal D$ is nonzero.
Therefore, $N(J+X)\ne0$.
\end{proof}

The following corollary is an application of \Cref{lem J-AB,prop UK}.

\begin{corollary}\label{cor UK}
Let $E,F,G$ be projections in $\mathcal M$, $A\in E\mathcal MF$, and $B\in G\mathcal MF$.
The following statements are equivalent:
\begin{enumerate}
\item[$(1)$] $R(A+YB)<E$ for every $Y\in E\mathcal MG$;
\item[$(2)$] $B=KA+L$ for some $K\in G\mathcal K_{\mathcal M}E$ and $L\in G\mathcal F_{\mathcal M}F$ with $R(L)\prec E-R(A)$.
\end{enumerate}
\end{corollary}

\begin{proof}
With the notation of \eqref{equ J-AB}, write $J=J_{A,B}$ and let $J=U|J|$ be the polar decomposition.
By \Cref{lem J-AB,prop UK}, (1) is equivalent to
\begin{enumerate}
\item[$(3)$] $N(J)\succ N(J^*)$ and $\chi_{(0,c]}(|J|)$ is finite for every $c>0$.
\end{enumerate}
Thus, it suffices to prove that $(2)\Leftrightarrow(3)$.
Recall that $N(J)=E-R(A)$.

$(3)\Rightarrow(2)$.
Since $N(J)\succ N(J^*)$, the projection $N(J^*)$ is finite.
Let
\begin{equation*}
	Q_c=N(J^*)+U\chi_{(0,c]}(|J|)U^*\quad (c>0).
\end{equation*}
Since finite projections $Q_c$ decrease strongly to $N(J^*)$ as $c\to0$, it is clear that $Q_\varepsilon\prec N(J)$ for sufficiently small $\varepsilon>0$.
Let $D=|J|^{-1}\chi_{(\varepsilon,\infty)}(|J|)$.
Since each $|J|^{-1}\chi_{(\varepsilon,n)}(|J|)$ is finite, $D$ is relatively compact.
Since
\begin{equation*}
	B^*(Q-Q_\varepsilon)
	=B^*JDU^*
	=A^*DU^*,
\end{equation*}
we have $B=UDA+Q_\varepsilon B$.
Thus, $K=UD\in G\mathcal K_{\mathcal M}E$ and $L=Q_\varepsilon B\in G\mathcal F_{\mathcal M}F$ satisfy
\begin{equation*}
	B=KA+L,\qquad
	R(L)\precsim Q_\varepsilon\prec N(J)=E-R(A).
\end{equation*}

$(2)\Rightarrow(3)$.
For every $\xi\in\ker J^*\cap\ker L^*$, we have $A^*K^*\xi=B^*\xi$.
It follows that $K^*\xi\in\mathcal D_{A,B}$ and $JK^*\xi=\xi$.
Since $\ker J^*\perp\operatorname{ran}J$, we have $\xi=0$.
Therefore, $L^*$ is injective on $\ker J^*$, and
\begin{equation*}
	N(J^*)\precsim R(L^*)\sim R(L)\prec E-R(A)=N(J).
\end{equation*}
Fix $c>0$ and let
\begin{equation*}
	G_c=\chi_{(0,c]}(|J|),\qquad
	S=L^*JG_c,\qquad
	H=G_c-R(S^*).
\end{equation*}
Since $JG_c$ is bounded and $L$ is finite, $S$ is finite and $R(S^*)\leqslant G_c$.
Moreover, it follows from $L^*JH=0$ that $A^*H=B^*JH=A^*K^*JH$.
Since $H\leqslant G_c\leqslant R(A)$, we have
\begin{equation*}
	H=R(A)K^*JH\in\mathcal K_{\mathcal M}.
\end{equation*}
Thus, $H$ is a relatively compact projection and hence finite.
Since $R(S^*)$ is also finite, $G_c=H+R(S^*)$ is finite.
This proves (3).
\end{proof}

We now obtain the following characterization.

\begin{theorem}\label{thm SR}
Let $\mathcal M$ be a separable factor and $T$ an operator in $\mathcal M$ with \eqref{equ assumption} and \eqref{equ T-matrix}.
Let $\mathcal D_T=\{\xi\in\operatorname{ran}E_T\colon A_T^*\xi\in\operatorname{ran}B_T^*\}$ and
\begin{equation*}
	J_T=(B_T^*|_{\operatorname{ran}Q})^{-1}A_T^*\colon\mathcal D_T(\subseteq\operatorname{ran}P)\to\operatorname{ran}Q,
\end{equation*}
where $P$ is the projection onto $\overline{\mathcal D_T}$ and $Q=R(B_T)$.
The following statements are equivalent:
\begin{enumerate}
\item[$(1)$] $T$ is strongly reducible in $\mathcal M$;
\item[$(2)$] $R(A_T+YB_T)<E_T$ for every $Y\in E_T\mathcal MG_T$;
\item[$(3)$] $B_T=KA_T+L$ for some relatively compact operator $K\in G_T\mathcal ME_T$ and finite operator $L\in G_T\mathcal MF_T$ satisfying that $R(L)\prec E_T-R(A_T)$;
\item[$(4)$] $N(J_T)\succ N(J_T^*)$ and $\chi_{(0,c]}(|J_T|)$ is finite for every $c>0$.
\end{enumerate}
\end{theorem}

\begin{proof}
The equivalence $(1)\Leftrightarrow(2)$ follows from \Cref{prop SR-criterion}.
Applying \Cref{cor UK} to $A=A_T$ and $B=B_T$ gives $(2)\Leftrightarrow(3)$.
Finally, $(2)\Leftrightarrow(4)$ follows from \Cref{lem J-AB,prop UK}.
\end{proof}

\begin{remark}\label{rem SR}
If $T$ is strongly reducible, then $B_T$ is relatively compact by (3).
\end{remark}

We present a concrete example to illustrate how to use \Cref{thm SR}.

\begin{example}\label{eg shift}
Let $\mathcal H$ be an infinite-dimensional separable Hilbert space and $T$ an operator in $\mathcal B(\mathcal H^{(4)})$ of the form
\begin{equation*}
	T=
	\left(\begin{array}{cc|c|c}
		0 & 0 & 0 & A\\
		\hline
		0 & 0 & C & D\\
		\hline
		0 & 0 & 0 & B\\
		0 & 0 & 0 & 0
	\end{array}\right),
\end{equation*}
where $A,B,C,D\in\mathcal B(\mathcal H)$ and $\ker B=\ker C=\ker C^*=\{0\}$.
Then $T^2\ne 0$, $T^3=0$, and
\begin{equation*}
	A_T=A,\quad B_T=
	\begin{pmatrix}
		B\\0
	\end{pmatrix},\quad C_T=C,\quad D_T=D.
\end{equation*}
Suppose that $\{e_n\}_{n=0}^{\infty}$ is an orthonormal basis for $\mathcal H$, $k_1,k_2\geqslant 0$ are integers, and $A,B$ are given by
\begin{equation*}
	Ae_n=a_ne_{n+k_1},\quad Be_n=b_ne_{n+k_2}\quad\text{for}~n\geqslant 0,
\end{equation*}
where $a_n\ne 0$ and $b_n\ne 0$.
Then $T$ is strongly reducible if and only if
\begin{equation}\label{equ k1-k2-a-b}
	k_1\geqslant 1\quad\text{and}\quad\lim_{n\to\infty}\frac{b_n}{a_n}=0.
\end{equation}
Indeed, with the notation in \Cref{thm SR}, we have
\begin{equation*}
	\operatorname{ran}Q=\overline{\mathrm{span}}\{e_{n+k_2}\colon n\geqslant 0\},\quad\mathcal D_T=\left\{\sum_{n=0}^{\infty}x_ne_n\colon\sum_{n=0}^{\infty}\frac{|a_n|^2}{|b_n|^2}|x_{n+k_1}|^2<\infty\right\},
\end{equation*}
and $P=I_{\mathcal H}$.
Moreover, $J_T$ is given by
\begin{equation*}
	J_T\colon\mathcal D_T\to\operatorname{ran}Q,\quad J_Te_n=0~(0\leqslant n<k_1),\quad J_Te_{n+k_1}=\frac{\overline{a_n}}{\overline{b_n}}e_{n+k_2}~(n\geqslant 0).
\end{equation*}
Therefore, condition $(4)$ in \Cref{thm SR} becomes \eqref{equ k1-k2-a-b}.
\end{example}

\section{Applications}
\label{sec application}

In this section, we present several applications of \Cref{thm SR}.
By \Cref{rem SR}, it suffices to consider the case that $B_T$ is a relatively compact operator in $\mathcal M$.
The following criterion applies when $A_T$ has closed range.

\begin{theorem}\label{thm SR-closed-range}
Let $\mathcal M$ be a separable semifinite factor and $T$ an operator in $\mathcal M$ with \eqref{equ assumption} and \eqref{equ T-matrix}.
If $B_T$ is relatively compact in $\mathcal M$ and $A_T$ has closed range, then the following statements are equivalent:
\begin{enumerate}
\item[$(1)$] $T$ is strongly reducible in $\mathcal M$;
\item[$(2)$] $N(A_T)\wedge F_T\prec N(A_T^*)\wedge E_T$.
\end{enumerate}
\end{theorem}

\begin{proof}
Let $F_0=N(A_T)\wedge F_T=F_T-R(A_T^*)$ and $E_0=N(A_T^*)\wedge E_T=E_T-R(A_T)$.
Since $R(B_T^*)=F_T$, it follows that $B_T$ is injective on $\operatorname{ran}F_T$.

$(1)\Rightarrow(2)$.
By \Cref{thm SR}, we have $B_T=KA_T+L$ for some relatively compact operator $K$ and finite operator $L$ with $R(L)\prec E_0$.
Since $A_TF_0=0$, it follows that $B_TF_0=LF_0$.
Thus, $F_0\sim R(B_TF_0)=R(LF_0)\precsim R(L)\prec E_0$.

$(2)\Rightarrow(1)$.
Since $F_0\prec E_0$, the projection $F_0$ is finite.
Since $A_T$ has closed range, there exists $X\in F_T\mathcal ME_T$ such that $XA_T=F_T-F_0$.
Let $K=B_TX$ and $L=B_TF_0$.
Then $K$ is relatively compact, $L$ is finite, and
\begin{equation*}
	B_T=KA_T+L,\qquad R(L)\sim F_0\prec E_0.
\end{equation*}
Therefore, $T$ is strongly reducible by \Cref{thm SR}.
\end{proof}

\begin{remark}\label{rem semi-Fredholm}
When $\mathcal M=\mathcal B(\mathcal H)$, condition $(2)$ means that $A_T$ is a semi-Fredholm operator with $\operatorname{ind}A_T<0$ (allowing $-\infty$) as an operator in $\mathcal B(F_T\mathcal H,E_T\mathcal H)$.
\end{remark}

If $A_T$ does not have closed range, it seems difficult to verify condition $(2)$ in \Cref{thm SR}.
In the following theorem, we consider the special case that $B_T$ is a finite operator, i.e., $F_T$ is a finite projection.
Note that
\begin{equation*}
	R(B_T)\sim R(B_T^*)=F_T=I-N(T^2)=R((T^*)^2)\sim R(T^2).
\end{equation*}

\begin{theorem}\label{thm SR-B-finite}
Let $\mathcal M$ be a separable semifinite factor and $T$ an operator in $\mathcal M$ satisfying condition \eqref{equ assumption}.
If $T^2$ is a finite operator, then the following statements are equivalent:
\begin{enumerate}
\item[$(1)$] $T$ is strongly reducible in $\mathcal M$;
\item[$(2)$] $R(T^2)\prec N(T)-R(TN(T^2))$; equivalently, $F_T\prec E_T$.
\end{enumerate}
\end{theorem}

\begin{proof}
Let $F_0=F_T-R(A_T^*)$ and $E_0=E_T-R(A_T)$.
Since $R(B_T^*)=F_T\sim R(T^2)$, $F_T$ is finite and $B_T$ is injective on $\operatorname{ran}F_T$.

$(1)\Rightarrow(2)$.
As in the proof of $(1)\Rightarrow(2)$ in \Cref{thm SR-closed-range}, we have $F_0\prec E_0$.
Since
\begin{equation*}
	F_T-F_0=R(A_T^*)\sim R(A_T)=E_T-E_0
\end{equation*}
and these projections are finite, it follows that $F_T\prec E_T$.

$(2)\Rightarrow(1)$.
For every $Y\in E_T\mathcal MG_T$, the operator $A_T+YB_T\in E_T\mathcal MF_T$ satisfies
\begin{equation*}
	R(A_T+YB_T)\precsim F_T\prec E_T.
\end{equation*}
Thus, $R(A_T+YB_T)<E_T$, and hence $T$ is strongly reducible by \Cref{thm SR}.
\end{proof}

Since $R(T)=I$ implies $N(T)\wedge N(T^*)=0$, the following corollary shows that a strongly reducible operator satisfying \eqref{equ assumption} cannot have dense range.

\begin{corollary}\label{cor N-T*}
Suppose that $T$ is an operator in $\mathcal M$ such that
\begin{equation*}
	N(T)\wedge N(T^*)=0,\quad T^2\ne 0,\quad\text{and}\quad\bigvee_{n=1}^{\infty}\ker T^n=\mathcal H.
\end{equation*}
Then $T$ is not strongly reducible.
\end{corollary}

\begin{proof}
Since $N(J_T)=E_T-R(A_T)=N(T)\wedge N(T^*)$ by \Cref{lem common-kernel}, we cannot have $N(J_T)\succ N(J_T^*)$.
Therefore, $T$ is not strongly reducible by \Cref{thm SR}.
\end{proof}

C. Jiang and J. Li proved in \cite[Theorem 4.1]{JL-00} that every Cowen-Douglas operator with finite index is not strongly reducible.
We show that the finite index assumption can be removed as follows.

\begin{corollary}\label{cor Cowen-Douglas}
Let $T$ be a Cowen-Douglas operator of index $1\leqslant n\leqslant\infty$.
Then $T$ is not strongly reducible.
\end{corollary}

\begin{proof}
Let $T\in B_n(\Omega)$, where $\Omega$ is a bounded connected open subset of $\mathbb{C}$.
After shifting an eigenvalue to zero, we may assume that $0\in\Omega$.
Then
\begin{equation*}
	\operatorname{ran}T=\mathcal H,\quad T^2\ne 0,\quad\text{and}\quad\bigvee_{n=1}^{\infty}\ker T^n=\mathcal H.
\end{equation*}
We complete the proof by \Cref{cor N-T*}.
\end{proof}

As usual, we denote by $J_n$ the canonical nilpotent Jordan block in $M_n(\mathbb{C})$, i.e.,
\begin{equation*}
	J_n=
	\begin{pmatrix}
		0 & 1 & & &\\
		& 0 & 1 & &\\
		& & \ddots & \ddots &\\
		& & & 0 & 1\\
		& & & & 0
	\end{pmatrix}_{n\times n}.
\end{equation*}
Let $S$ be the unilateral shift operator in $\mathcal B(\ell^2)$.
Since $S^*$ has a similar (infinite) matrix representation, we write $J_\infty=S^*$ for simplicity.

\begin{corollary}\label{cor Jordan}
Suppose that
\begin{equation*}
	T=0^{(n_1)}\oplus J_2^{(n_2)}\oplus J_3^{(n_3)}\oplus\cdots\oplus J_\infty^{(n_\infty)}
\end{equation*}
is an operator acting on a separable Hilbert space $\mathcal H$ with $T^2\ne 0$, where $0\leqslant n_j\leqslant\infty$ for all $1\leqslant j\leqslant\infty$.
The following statements are equivalent:
\begin{enumerate}
\item[$(1)$] $T$ is strongly reducible in $\mathcal B(\mathcal H)$;
\item[$(2)$] $n_1>\sum_{j=3}^{\infty}(j-2)n_j$ and $n_\infty=0$.
\end{enumerate}
\end{corollary}

\begin{proof}
A direct computation shows that
\begin{equation*}
	\operatorname{rank}E_T=n_1\quad\text{and}\quad\operatorname{rank}F_T=\sum_{j=3}^{\infty}(j-2)n_j+\infty\cdot n_\infty.
\end{equation*}
Since $A_T=0$, it has closed range.
If $T$ is strongly reducible, then $B_T$ is compact by \Cref{rem SR} and $F_T\prec E_T$ by \Cref{thm SR-closed-range}.
It follows that $n_1>\sum_{j=3}^{\infty}(j-2)n_j$ and $n_\infty=0$.
Conversely, if condition $(2)$ holds, then $B_T$ is finite-rank and $F_T\prec E_T$.
Thus, $T$ is strongly reducible by \Cref{thm SR-closed-range}.
This completes the proof.
\end{proof}

\begin{remark}\label{rem Jordan}
Note that strong reducibility is invariant under taking adjoints.
As a consequence of \Cref{cor Jordan}, $0^{(n)}\oplus S$ is not strongly reducible for every $0\leqslant n\leqslant\infty$, where $S$ is the unilateral shift operator.
\end{remark}

\section{Algebraic operators}
\label{sec algebric}

In this section, we study the strong reducibility of algebraic operators in finite factors and type $\mathrm{III}$ factors.
Recall that an operator $T$ is called \emph{algebraic} if $p(T)=0$ for some nonzero polynomial $p$.

\subsection{Finite factors}
\label{subsec finite}

Recall that factors of type $\mathrm{I}_n$ and type $\mathrm{II}_1$ are called \emph{finite factors}, and there exists a unique normal faithful tracial state on every finite factor.
We have the following simple trace criterion for the strong reducibility of operators satisfying \eqref{equ assumption}.

\begin{corollary}\label{cor SR-finite-factor}
Let $(\mathcal M,\tau)$ be a separable finite factor and $T$ an operator in $\mathcal M$ satisfying condition \eqref{equ assumption}.
Then $T$ is strongly reducible if and only if $\tau(R(T))<\frac{1}{2}$.
\end{corollary}

\begin{proof}
By \Cref{thm SR-B-finite}, $T$ is strongly reducible if and only if $F_T\prec E_T$.
With the notation in \eqref{equ kernel-flag}, we have $\tau(F_T)=1-\tau(P_1)-\tau(P_2)$ and
\begin{equation*}
	\tau(E_T)=\tau(P_1)-\tau(R(T_{12}))=\tau(P_1)-\tau(R(T_{12}^*))=\tau(P_1)-\tau(P_2).
\end{equation*}
Thus, $F_T\prec E_T$ if and only if $\tau(P_1)>\frac{1}{2}$, i.e., $\tau(R(T))<\frac{1}{2}$.
\end{proof}

Let $P$ be a projection in a factor $\mathcal M$.
For every operator $A\in P\mathcal MP$, its kernel in $\operatorname{ran}P$ is denoted by
\begin{equation*}
	N_P(A)=N(A)\wedge P.
\end{equation*}
In the following lemma, we consider the direct sum of two nilpotent operators.

\begin{lemma}\label{lem A-V-B}
Let $(\mathcal M,\tau)$ be a finite factor and $P,Q$ nonzero projections in $\mathcal M$ with $PQ=0$.
Then for any nilpotent operators $A\in P\mathcal MP$ and $B\in Q\mathcal MQ$, there exists an operator $V\in P\mathcal MQ$ such that
\begin{equation*}
	C=A+V+B=
	\begin{pmatrix}
		A & V\\
		0 & B
	\end{pmatrix}
\end{equation*}
is nilpotent and $\tau(N_{P+Q}(C))=\max\{\tau(N_P(A)),\tau(N_Q(B))\}$.
Moreover, if $A\ne 0$, then $V$ can be chosen such that $C^2\ne 0$.
\end{lemma}

\begin{proof}
Let $a=\tau(N_P(A))=\tau(N_P(A^*))$ and $b=\tau(N_Q(B))$.
Let $V$ be a partial isometry in $P\mathcal MQ$ such that
\begin{equation*}
	P_0=VV^*\leqslant N_P(A^*),\qquad
	Q_0=V^*V\leqslant N_Q(B),\qquad
	\tau(P_0)=\tau(Q_0)=\min\{a,b\}.
\end{equation*}
If $A^r=0$ and $B^s=0$, then $C^{r+s}=0$.
Since $A^*V=0$ and $Q_0\leqslant N_Q(B)$, we have
\begin{equation*}
	N_{P+Q}(C)=N_P(A)+N_Q(B)-Q_0.
\end{equation*}
In particular, $\tau(N_{P+Q}(C))=\max\{a,b\}$.
Now suppose that $A\ne 0$.
If $A^2\ne 0$, then clearly $C^2\ne 0$.
If $A^2=0$, then $0\ne R(A^*)\leqslant N_P(A^*)$.
In this case, we can choose $V$ such that $AP_0\ne 0$, i.e., $AV\ne 0$.
Since $BQ_0=0$ and $VQ_0=V$, it follows that $C^2Q_0=AV\ne0$.
This proves the last assertion.
\end{proof}

The following technical lemma will be used in the proof of \Cref{thm Hsin}.

\begin{lemma}\label{lem D+N}
Let $(\mathcal M,\tau)$ be a finite factor, $\{P_j\}_{j=1}^n$ nonzero projections in $\mathcal M$ with sum $I$, $\{A_j\}_{j=1}^{n}$ nilpotent operators with $A_j\in P_j\mathcal MP_j$, $\{\lambda_j\}_{j=1}^n$ distinct scalars, and
\begin{equation*}
	T=\sum_{j=1}^n(\lambda_jP_j+A_j).
\end{equation*}
Suppose that, with respect to the decomposition $I=\sum_{j=1}^nP_j$, there exists an upper triangular operator 
\begin{equation*}
	A=
	\begin{pmatrix}
		A_{11} & A_{12} & \cdots & A_{1n}\\
		& A_{22} & \cdots & A_{2n}\\
		& & \ddots & \vdots\\
		& & & A_{nn}
	\end{pmatrix}
\end{equation*}
such that $A$ is not strongly reducible and $A_{jj}=A_j$ for all $1\leqslant j\leqslant n$.
Then $T$ is not strongly reducible.
\end{lemma}

\begin{proof}
Let $X$ be an invertible operator in $\mathcal M$ such that $XAX^{-1}$ is irreducible.
Let
\begin{equation*}
	E_j=P_1+\cdots+P_j\qquad\text{and}\qquad F_j=R(XE_j)\qquad\text{for}~1\leqslant j\leqslant n.
\end{equation*}
It is clear that $\tau(E_j)=\tau(F_j)$.
Thus, there exists a unitary operator $U$ in $\mathcal M$ such that $UE_jU^*=F_j$ for every $1\leqslant j\leqslant n$.
Let $Y=U^*X$.
Since $Y\operatorname{ran}E_j=\operatorname{ran}E_j$, both $Y$ and $Y^{-1}$ are upper triangular with respect to the decomposition $I=\sum_{j=1}^nP_j$.
Moreover, $YAY^{-1}=U^*(XAX^{-1})U$ is irreducible.
We denote by $\{Y_j\}_{j=1}^n$ the diagonal elements of $Y$.
Let
\begin{equation*}
	T_1=\sum_{j=1}^n(\lambda_jP_j+Y_jA_jY_j^{-1})\quad\text{and}\quad T_2=YAY^{-1}+\sum_{j=1}^n\lambda_jP_j.
\end{equation*}
Then $T_2$ is upper triangular with diagonal elements $\{\lambda_jP_j+Y_jA_jY_j^{-1}\}_{j=1}^n$.
Clearly, $T$ is similar to $T_1$.
Since the diagonal blocks $\{\lambda_jP_j+Y_jA_jY_j^{-1}\}_{j=1}^n$ have pairwise disjoint spectra, $T_1$ and $T_2$ are similar in $\mathcal M$ by \cite{Rosenblum-56}.

It suffices to prove that $T_2$ is irreducible in $\mathcal M$.
Suppose that $A_j^{m_j}=0$ for some integer $m_j\geqslant 1$.
Then
\begin{equation*}
	E_j=N\left(\prod_{k=1}^j(T_2-\lambda_k I)^{m_k}\right)\qquad\text{for every}~1\leqslant j\leqslant n.
\end{equation*}
Let $Q$ be a projection in $\mathcal M$ commuting with $T_2$.
Then $Q$ commutes with each $E_j$, hence with every $P_j$ and with $YAY^{-1}$.
It follows that $Q=0$ or $I$ by the irreducibility of $YAY^{-1}$.
Thus, $T_2$ is irreducible.
This completes the proof.
\end{proof}

In the following theorem, we give a complete criterion for the strong reducibility of algebraic operators in finite factors.

\begin{theorem}\label{thm Hsin}
Let $(\mathcal M,\tau)$ be a separable finite factor with $\dim\mathcal M>4$ and $T$ an algebraic operator in $\mathcal M$.
Then $T$ is not strongly reducible if and only if
\begin{enumerate}
\item[$(1)$] $\tau(R(\lambda I-T))\geqslant\frac{1}{2}$ for every scalar $\lambda$ in $\mathbb{C}$;
\item[$(2)$] the set $\{I,T,T^2\}$ is linearly independent over $\mathbb{C}$.
\end{enumerate}
\end{theorem}

\begin{proof}
The necessity follows from \cite[Lemma 4.3]{MSY-26}.
Suppose that conditions $(1)$ and $(2)$ hold.
Since $T$ is algebraic, up to similarity, we may assume that
\begin{equation*}
	T=\sum_{j=1}^n(\lambda_jP_j+A_j),
\end{equation*}
where $\{P_j\}_{j=1}^n$ are nonzero projections with sum $I$, $\{\lambda_j\}_{j=1}^n$ are distinct scalars, and each $A_j\in P_j\mathcal MP_j$ is nilpotent.
Since similarity preserves kernel traces, condition $(1)$ gives
\begin{equation*}
	\tau(N_{P_j}(A_j))=\tau(N(T-\lambda_jI))\leqslant\frac{1}{2}
	\qquad (1\leqslant j\leqslant n).
\end{equation*}
If $n=1$, then $T-\lambda_1I$ is nilpotent of order at least $3$ by condition $(2)$, and the conclusion follows from \Cref{cor SR-finite-factor}.
If $A_j=0$ for every $1\leqslant j\leqslant n$, then $T$ is normal, and the conclusion follows from \cite[Theorem 1.2]{MSY-26}.

Suppose that $n\geqslant 2$ and without loss of generality that $A_1\ne 0$.
Repeated application of \Cref{lem A-V-B} produces a nilpotent operator $A\in\mathcal M$, upper triangular with respect to $I=\sum_{j=1}^nP_j$, with diagonal elements $\{A_j\}_{j=1}^n$, such that $A^2\ne0$ and
\begin{equation*}
	\tau(N(A))=\max_{1\leqslant j\leqslant n}\tau(N_{P_j}(A_j))\leqslant\frac{1}{2}.
\end{equation*}
Thus, $A$ is not strongly reducible by \Cref{cor SR-finite-factor}, and the conclusion follows from \Cref{lem D+N}.
\end{proof}

\subsection{Type \texorpdfstring{$\mathrm{III}$}{III} factors}
\label{subsec III}

Since every nonzero projection in a type $\mathrm{III}$ factor is infinite, we have $\mathcal K_{\mathcal M}=\{0\}$.
Hence we can immediately get the following result by \Cref{thm SR}.

\begin{corollary}\label{cor SR-III}
Let $\mathcal M$ be a separable type $\mathrm{III}$ factor and $T$ an operator in $\mathcal M$ satisfying condition \eqref{equ assumption}.
Then $T$ is not strongly reducible.
\end{corollary}

\begin{proof}
If $T$ were strongly reducible, then condition $(3)$ in \Cref{thm SR} would give $B_T=0$.
However, $R(B_T^*)=F_T=I-N(T^2)\ne0$ as $T^2\ne0$.
This contradiction completes the proof.
\end{proof}

For algebraic operators in type $\mathrm{III}$ factors, the criterion in \Cref{thm Hsin} reduces to the absence of a quadratic equation.

\begin{theorem}\label{thm SR-III}
Let $\mathcal M$ be a separable type $\mathrm{III}$ factor and $T$ an algebraic operator in $\mathcal M$.
Then $T$ is not strongly reducible if and only if the set $\{I,T,T^2\}$ is linearly independent over $\mathbb{C}$.
\end{theorem}

\begin{proof}
Since the proof is analogous to that of \Cref{thm Hsin}, we will only outline the strategy of the proof.
The necessity follows from \cite[Lemma 4.3]{MSY-26}.
For the sufficiency, suppose that $\{I,T,T^2\}$ is linearly independent over $\mathbb C$.
First, \Cref{lem A-V-B} remains valid in type $\mathrm{III}$ factors with the trace condition omitted, including its last assertion concerning $C^2\ne0$.
Moreover, \Cref{lem D+N} also holds in type $\mathrm{III}$ factors.
Actually, setting $F_0=0$ in its proof, the projections $P_j$ and $F_j-F_{j-1}$ are nonzero and hence equivalent for every $1\leqslant j\leqslant n$.
Thus, there exists a unitary operator $U$ in $\mathcal M$ such that $UE_jU^*=F_j$ for every $1\leqslant j\leqslant n$.
The remainder of that proof applies unchanged.
Therefore, by applying \Cref{cor SR-III} and the proof of \Cref{thm Hsin}, we can get the desired conclusion.
\end{proof}

\section{Multi-cyclic operators}
\label{sec multi-cyclic}

In this section, we prove that every operator with an infinite-dimensional multi-cyclic direct summand is similar to an irreducible operator.
Recall that $S\in\mathcal{B}(\mathcal{H})$ is multi-cyclic if there is a generating set $\{\xi_j\}_{j=1}^m$ in $\mathcal H$ such that
\begin{equation*}
	\mathcal{H}=\overline{\mathrm{span}}\{S^k\xi_j\colon k\geqslant 0,~1\leqslant j\leqslant m\}.
\end{equation*}
We first show that the generating set can be chosen such that their first few iterates by $S$ are linearly independent.

\begin{lemma}\label{lem independent-generators}
Let $S$ be a multi-cyclic operator acting on an infinite-dimensional Hilbert space $\mathcal{H}$.
Then $S$ has a generating set $\{\xi_j\}_{j=1}^m$ such that the set
\begin{equation}\label{equ S^k-xi}
	\{S^k\xi_j\colon 0\leqslant k\leqslant 3,~1\leqslant j\leqslant m\}
\end{equation}
is linearly independent.
\end{lemma}

\begin{proof}
Let $\{\xi_j\}_{j=1}^m$ be a generating set for $S$.
Since $\mathcal{H}$ is infinite dimensional, after reordering, we may assume that the orbit $\{S^n\xi_1\}_{n=0}^\infty$ of $\xi_1$ is linearly independent.
Let 
\begin{equation*}
	\mathcal{H}_0=\mathrm{span}\{S^k\xi_j\colon 0\leqslant k\leqslant 3,~1\leqslant j\leqslant m\}.
\end{equation*}
Since $\mathcal{H}_0$ is finite dimensional, there exists $d\geqslant 3$ such that
\begin{equation*}
	\mathcal{H}_0\cap\mathrm{span}\{S^n\xi_1\colon n\geqslant 0\}\subseteq\mathrm{span}\{S^n\xi_1\colon 0\leqslant n\leqslant d\}.
\end{equation*}
Choose integers $N_2,\ldots,N_m>d$ such that the intervals $[N_j,N_j+3]$ are pairwise disjoint.
Then the vectors $\{S^{N_j+k}\xi_1\colon 0\leqslant k\leqslant 3,~2\leqslant j\leqslant m\}$ are linearly independent modulo $\mathcal{H}_0$.
Replacing $\xi_j$ by $\xi_j+S^{N_j}\xi_1$ for $2\leqslant j\leqslant m$ preserves the generating property, since $\xi_1$ is unchanged.
With this new generating set, the set \eqref{equ S^k-xi} is linearly independent.
\end{proof}

The following lemma provides a convenient generating set up to similarity.

\begin{lemma}\label{lem nice}
Let $S$ be a multi-cyclic operator acting on an infinite-dimensional Hilbert space $\mathcal{H}$.
Then there exists an operator $S_0$ similar to $S$ and a generating set $\{\xi_j\}_{j=1}^m$ for $S_0$ such that
\begin{enumerate}
\item[$(1)$] $\{\xi_j\}_{j=1}^m\cup\{S_0\xi_j\}_{j=1}^m$ is orthonormal and $\langle S_0^2\xi_i,S_0^2\xi_j\rangle\ne 0$ for all $i\ne j$;
\item[$(2)$] $S_0^*S_0\xi_j=\xi_j$ for $1\leqslant j\leqslant m$;
\item[$(3)$] there exists an orthogonal decomposition $\mathcal{H}=\mathcal{H}_1\oplus\mathcal{H}_2\oplus\mathcal{H}_3$ such that
\begin{equation*}
	\mathcal{H}_1=\mathrm{span}\{\xi_j\}_{j=1}^m,\qquad\mathcal{H}_2=\mathrm{span}\{S_0\xi_j\}_{j=1}^m
\end{equation*}
and
\begin{equation*}
	S_0\mathcal{H}_1=\mathcal{H}_2,\qquad S_0\mathcal{H}_2\subseteq\mathcal{H}_3,\qquad S_0\mathcal{H}_3\subseteq\mathcal{H}_1\oplus\mathcal{H}_3.
\end{equation*}
\end{enumerate}
\end{lemma}

\begin{proof}
Let $\{\xi_j\}_{j=1}^m$ be a generating set for $S$ given by \Cref{lem independent-generators}.
Let $\{e_j\}_{j=1}^m$ be an orthonormal basis for $\mathbb C^m$ and choose a bounded linear operator $\Phi\colon\mathcal H\to\mathbb C^m$ such that
\begin{equation*}
	\Phi\xi_j=\Phi S^2\xi_j=\Phi S^3\xi_j=0,\qquad
	\Phi S\xi_j=e_j\quad(1\leqslant j\leqslant m).
\end{equation*}
Let $\mathcal{H}_1=\mathrm{span}\{\xi_j\}_{j=1}^m$, $\mathcal{H}_2=\mathrm{span}\{S\xi_j\}_{j=1}^m$, and $\mathcal{H}_3=\ker\Phi\cap\ker\Phi S$.
Then
\begin{equation*}
	\mathcal H=\mathcal H_1\dotplus\mathcal H_2\dotplus\mathcal H_3,\qquad
	S\mathcal H_1=\mathcal H_2,\qquad
	S\mathcal H_2\subseteq\mathcal H_3,\qquad
	S\mathcal H_3\subseteq\mathcal H_1+\mathcal H_3.
\end{equation*}
Since the set \eqref{equ S^k-xi} is linearly independent, there exists an equivalent Hilbert inner product $\langle\cdot,\cdot\rangle_{\mathrm{new}}$ on $\mathcal H$ such that the above three summands are orthogonal, the set $\{\xi_j\}_{j=1}^m\cup\{S\xi_j\}_{j=1}^m$ is orthonormal, and
\begin{equation*}
	\langle S^2\xi_i,S^2\xi_j\rangle_{\mathrm{new}}\ne0
	\qquad\text{for all}~i\ne j.
\end{equation*}
Let $A\in\mathcal B(\mathcal H)$ be a positive invertible operator such that
\begin{equation*}
	\langle x,y\rangle_{\mathrm{new}}=\langle Ax,Ay\rangle.
\end{equation*}
Let $S_0=ASA^{-1}$, $\xi_j'=A\xi_j$, and $\mathcal{H}_k'=A\mathcal{H}_k$.
It is clear that conditions $(1)$ and $(3)$ hold for the $\{\xi_j'\}_{j=1}^m$ and $\{\mathcal H_j'\}_{j=1}^3$.
Since $S_0(\mathcal H_2'\oplus\mathcal H_3')\subseteq\mathcal H_1'\oplus\mathcal H_3'$, we have
\begin{equation*}
	\langle S_0^*S_0\xi_j',x\rangle
	=\langle S_0\xi_j',S_0x\rangle
	=
	\begin{cases}
		\delta_{ij}, & x=\xi_i',\\
		0, & x\in\mathcal H_2'\oplus\mathcal H_3'.
	\end{cases}
\end{equation*}
Thus, $S_0^*S_0\xi_j'=\xi_j'$ for every $1\leqslant j\leqslant m$.
This completes the proof.
\end{proof}

The next lemma constructs an injective operator whose range is orthogonal to a prescribed finite-dimensional subspace.

\begin{lemma}\label{lem injective-coupling}
Let $\mathcal H$ be a separable infinite-dimensional Hilbert space, $\mathcal F\subseteq\mathcal H$ a finite-dimensional subspace, and $S\in\mathcal B(\mathcal H)$ not a scalar plus a finite-rank operator.
Then for every $T\in\mathcal B(\mathcal K)$ with $\mathcal K$ separable, there exists an isometry $V\colon\mathcal{K}\to\mathcal{H}$ such that $VT-SV\colon\mathcal{K}\to\mathcal{H}$ is injective and $\operatorname{ran}V+\operatorname{ran}SV\subseteq\mathcal{F}^\perp$.
\end{lemma}

\begin{proof}
For every finite-dimensional subspace $\mathcal E\subseteq\mathcal H$, there is a unit vector $u$ such that
\begin{equation*}
	u\perp\mathcal E+S^*\mathcal E,\qquad Su\notin\mathbb Cu.
\end{equation*}
Otherwise $S$ would act as a scalar $\lambda$ on $(\mathcal E+S^*\mathcal E)^\perp$.
Since this subspace has finite codimension, $S-\lambda I$ would have finite rank, a contradiction.
By taking $\mathcal E=\mathcal F$, we can choose a unit vector $u_1\perp\mathcal F+S^*\mathcal F$ such that
\begin{equation*}
	Su_1=\lambda_1u_1+\mu_1v_1
\end{equation*}
for some scalar $\mu_1>0$ and some unit vector $v_1\perp u_1$.
Then $\{u_1,v_1\}$ is an orthonormal set in $\mathcal F^\perp$.
By taking $\mathcal E=\mathrm{span}(\mathcal F\cup\{u_1,v_1\})$, we can similarly obtain an orthonormal set $\{u_2,v_2\}$ in $(\mathcal F\cup\{u_1,v_1\})^\perp$ such that $Su_2=\lambda_2u_2+\mu_2v_2$ with $\mu_2>0$.
Inductively, we obtain an orthonormal set $\{u_n,v_n\colon n\geqslant1\}$ in $\mathcal F^\perp$ such that
\begin{equation*}
	Su_n=\lambda_nu_n+\mu_nv_n,\qquad \mu_n>0.
\end{equation*}
Choose an isometry $V\colon\mathcal K\to\overline{\mathrm{span}}\{u_n\colon n\geqslant1\}$.
Then $\operatorname{ran}V+\operatorname{ran}SV\subseteq\mathcal F^\perp$.
For every $x\in\ker(VT-SV)$, we have
\begin{equation*}
	0=\langle(VT-SV)x,v_n\rangle
	=-\mu_n\langle Vx,u_n\rangle
	\qquad(n\geqslant1).
\end{equation*}
Since each $\mu_n>0$ and $Vx\in\overline{\mathrm{span}}\{u_1,u_2,u_3,\ldots\}$, it follows that $Vx=0$, and hence $x=0$.
Thus, $VT-SV$ is injective.
\end{proof}

We now give an irreducibility criterion for triangular operators with an injective off-diagonal entry.

\begin{lemma}\label{lem triangular-irreducibility}
Suppose that $W\in\mathcal B(\mathcal H)$ admits an orthonormal generating set $\{\xi_j\}_{j=1}^m$, $B\in\mathcal B(\mathcal K,\mathcal H)$ is injective, $T\in\mathcal{B}(\mathcal K)$, and
\begin{equation*}
	A=\begin{pmatrix}T&0\\B&W\end{pmatrix}\in\mathcal{B}(\mathcal K\oplus\mathcal H).
\end{equation*}
Let $e_j=(0,\xi_j)$ for $1\leqslant j\leqslant m$.
If there exist distinct scalars $\{\lambda_j\}_{j=1}^m$ such that
\begin{equation*}
	\ker(A^*A-\lambda_jI)=\mathbb Ce_j\quad(1\leqslant j\leqslant m),\qquad
	\langle W^2\xi_i,W^2\xi_j\rangle\ne0\quad (i\ne j),
\end{equation*}
then $A$ is irreducible.
\end{lemma}

\begin{proof}
Let $Q$ be a projection commuting with $A$.
Then $Q$ also commutes with $A^*$ and $A^*A$.
It follows that
\begin{equation*}
	Qe_j=\varepsilon_je_j,\qquad\text{where}~\varepsilon_j\in\{0,1\}.
\end{equation*}
Since $Q$ commutes with $A^2$ and $A^2e_j=(0,W^2\xi_j)$, we have
\begin{equation*}
	(\varepsilon_i-\varepsilon_j)\langle W^2\xi_i,W^2\xi_j\rangle=0\qquad\text{for all}~1\leqslant i,j\leqslant m.
\end{equation*}
Thus, all $\varepsilon_j$ are equal.
By considering $I-Q$ if necessary, we may assume that $Qe_j=0$ for every $j$.
Since $\{\xi_j\}_{j=1}^m$ generates $\mathcal{H}$ under $W$, it follows that $Q|_{\{0\}\oplus\mathcal{H}}=0$.
Hence $Q=Q_1\oplus0$ for a projection $Q_1$ on $\mathcal{K}$.
The identity $QA=AQ$ gives $BQ_1=0$, and hence $Q_1=0$ by the injectivity of $B$.
Therefore, $Q=0$ and $A$ is irreducible.
\end{proof}

We now combine the preceding lemmas to prove the following theorem.

\begin{theorem}\label{thm multi-cyclic}
Let $S$ be a multi-cyclic operator acting on an infinite-dimensional Hilbert space $\mathcal{H}$.
Then $T\oplus S$ is not strongly reducible for every operator $T$ acting on a separable Hilbert space.
In particular, $S$ is not strongly reducible.
\end{theorem}

\begin{proof}
Let $T\in\mathcal{B}(\mathcal{K})$ with $\mathcal{K}$ separable.
We may assume that $S$ has a generating set $\{\xi_j\}_{j=1}^m$ and orthogonal summands $\mathcal{H}_1$, $\mathcal{H}_2$, and $\mathcal{H}_3$ as in \Cref{lem nice}.
Observe that $S$ is not a scalar plus a finite-rank operator, since otherwise all the generating orbits would lie in a finite-dimensional space.
By \Cref{lem injective-coupling} applied with $\mathcal F=\mathcal H_1\oplus\mathcal H_2$, there exists an isometry $V\colon\mathcal{K}\to\mathcal{H}$ such that $B=VT-SV$ is injective and
\begin{equation}\label{equ ran-B}
	\operatorname{ran}B\subseteq\operatorname{ran}V+\operatorname{ran}SV\subseteq\mathcal{H}_3.
\end{equation}
Choose distinct positive scalars $a_1,\ldots,a_m>1+\|T\|+\|B\|+\|S\|$ and define an invertible operator $D\in\mathcal{B}(\mathcal{H})$ by
\begin{equation*}
	D\xi_j=a_j^{-1}\xi_j\quad(1\leqslant j\leqslant m),\qquad D|_{\mathcal{H}_2\oplus\mathcal{H}_3}=I_{\mathcal{H}_2\oplus\mathcal{H}_3}.
\end{equation*}
Let $W=DSD^{-1}$.
Since $D$ is a contraction and is the identity on $\operatorname{ran}V+\operatorname{ran}SV$, we have
\begin{equation}\label{equ W-norm}
	WV=SV,\qquad\|W|_{\mathcal{H}_2\oplus\mathcal{H}_3}\|\leqslant\|S\|.
\end{equation}
Let
\begin{equation*}
	X=\begin{pmatrix}I_{\mathcal{K}}&0\\V&D\end{pmatrix},\qquad A=X(T\oplus S)X^{-1}=\begin{pmatrix}T&0\\B&W\end{pmatrix}.
\end{equation*}
Let $e_j=(0,\xi_j)$ for $1\leqslant j\leqslant m$.
By the choice of $D$, \Cref{lem nice}, and \eqref{equ ran-B}, we have
\begin{equation*}
	A^*Ae_j=a_j^2e_j.
\end{equation*}
Moreover, \eqref{equ W-norm} gives $\|A|_{\mathcal{K}\oplus\mathcal{H}_2\oplus\mathcal{H}_3}\|\leqslant\|T\|+\|B\|+\|S\|<\min_{1\leqslant j\leqslant m}a_j$.
Since $\{a_j\}_{j=1}^m$ are distinct, it follows that
\begin{equation*}
	\ker(A^*A-a_j^2I)=\mathbb{C}e_j\qquad\text{for all}~1\leqslant j\leqslant m.
\end{equation*}
The set $\{\xi_j\}_{j=1}^m$ still generates $\mathcal{H}$ under $W$, since $D$ preserves $\mathcal{H}_1$.
Moreover, we have $W^2\xi_j=a_jS^2\xi_j$, and hence
\begin{equation*}
	\langle W^2\xi_i,W^2\xi_j\rangle=a_ia_j\langle S^2\xi_i,S^2\xi_j\rangle\ne0\qquad\text{for all}~i\ne j.
\end{equation*}
Therefore, $A$ is irreducible by \Cref{lem triangular-irreducibility}, and $T\oplus S$ is not strongly reducible.
\end{proof}

\begin{remark}\label{rem multicyclic}
Every injective unilateral weighted shift operator on $\ell^2$ is cyclic and the Volterra operator on $L^2[0,1]$ is also cyclic.
An operator on an infinite-dimensional separable Hilbert space with dense range and dense generalized kernel is also cyclic by \cite[Corollary 3.3]{BBP-04}.
Consequently, \Cref{thm multi-cyclic} applies whenever one of these operators occurs as a direct summand.
\end{remark}

\section*{Data Availability}
Data sharing is not applicable to this article, as no new data were created or analyzed in this work.

\section*{Ethics Declarations}
\noindent{\bf Conflict of interest}
On behalf of all authors, the corresponding author states that there is no conflict of interest.

\section*{AI Declarations}
The authors declare that no Generative AI or AI-assisted technologies were used in the writing, proofreading, or mathematical derivations of this manuscript.

\section*{Acknowledgments}
The authors gratefully acknowledge the anonymous referees for their insightful comments and constructive suggestions. 

\end{document}